\documentclass[11pt,a4paper]{article}

\pdfoutput=1

\usepackage[utf8]{inputenc}
\usepackage{authblk}
\usepackage{graphicx}
\usepackage{amsmath}
\usepackage{amssymb}
\usepackage{caption}
\usepackage{subcaption}
\usepackage{xcolor}
\usepackage{hyperref}
\usepackage{fancyhdr}

\hypersetup{
  colorlinks,
  linkcolor={red!70!black},
  citecolor={blue!50!black},
  urlcolor={blue!80!black}
}

\author[1,2]{Dimitris Vartziotis}
\author[ ]{Joachim Wipper}
\affil[1]{\small TWT GmbH Science \& Innovation, Department for Mathematical
  Research, Ernsthaldenstraße~17, 70565~Stuttgart, Germany}
\affil[2]{\small NIKI Ltd.\ Digital Engineering, Research Center, Ethnikis
  Antistasis~205, 45500~Katsikas, Ioannina, Greece}

\begin{document}
\title{The fractal nature of an approximate\\ prime counting function}
\date{November~7, 2016}
\maketitle

\begin{abstract}
  Prime number related fractal polygons and curves are derived by combining
  two different aspects. One is an approximation of the prime counting
  function build on an additive function. The other are prime number indexed
  basis entities taken from the discrete or continuous Fourier basis.
\end{abstract}

\section{Introduction}

The distribution of prime numbers is one of the central problems in analytic
number theory. Here, the prime counting function $\pi(x)$, giving the
number of primes less or equal to $x\in\mathbb{R}$, is of special interest.  As
stated in \cite{Ore1988}: {\em when the large-scale distribution of primes is
  considered, it appears in many ways quite regular and obeys simple
  laws}. One of the first central results regarding the asymptotic
distribution of primes is given by the prime number theorem, providing the limit
\begin{equation}
  \label{eq:piapproxbyln}
  \lim_{x\to\infty}\frac{\pi(x)}{x/\ln x} = 1,  
\end{equation}
which was proved independently in 1896 by Jacques Salomon Hadamard and
Charles-Jean de La Vall{\'e}e Poussin. Both proofs are based on complex
analysis using the Riemann zeta function $\zeta(s):=\sum_{n=1}^{\infty} 1/n^s$,
with $s\in \mathbb{C}$ and the fact that $\zeta(s)\neq 0$ for all
$s:=1+\mathrm{i}y$, $y>0$.

An improved approximation of $\pi(x)$ is given by the Eulerian logarithmic
integral
\begin{equation}
  \label{eq:Li}
  \operatorname{Li}(x):=\int_{2}^x\frac{{\rm d}t}{\ln t}\,.  
\end{equation}
This result was first mentioned by Carl Friedrich Gauß in 1849 in a letter to
Encke refining the estimate $n/\ln(n)$ of $\pi(n)$ given by the only 15 years
old Gauß in 1792. This conjecture was also stated by Legendre in 1798.

In 2007 Terence Tao gave an informal sketch of proof in his lecture
``Structure and randomness in the prime numbers''  as follows \cite{Tao2007}:
\begin{itemize}
\item {\em Create a ``sound wave'' (or more precisely, the von Mangoldt
    function) which is noisy at prime number times, and quite at other times.}
  [...]
\item {\em ``Listen'' (or take Fourier transforms) to this wave and record the
    notes that you hear (the zeroes of the Riemann zeta function, or the
    ``music of the primes''). Each such note corresponds to a hidden pattern
    in the distribution of the primes.}
\end{itemize}
In the same spirit, the present work tries to paint a picture of the
primes. By combining an alternative approximation of the prime counting
function $\pi(x)$ based on an additive function as proposed in
\cite{VartziotisTzavellas2016} with prime number related Fourier polygons used
in the context of regularizing polygon transformations as given in
\cite{VartziotisWipper2009PP, VartziotisWipper2010LT}, fractal prime polygons
and fractal prime curves are derived.

Three types of structures in the distribution of prime numbers are distinguished
in \cite{Batchko2014}. The first is {\em local structure}, like residue
classes or arithmetic progression \cite{GreenTao2008}. The second is {\em
  asymptotic structure} as provided by the prime number theorem. The third is
{\em statistical structure} as described for example in \cite{Cattani2010}
reporting an empirical evidence of fractal behavior in the distribution of
primes or \cite{Selvam2014} describing fractal fluctuations in the spacing
intervals of adjacent prime numbers generic to diverse dynamical systems in
nature. Quasi self similar structures in the distribution of differences of
prime-indexed primes with scaling by prime-index order have been observed in
\cite{Batchko2014}. In the present work, the asymptotic structure becomes
visually apparent by the given fractals.

\section{Approximations of the prime counting function}
\label{sec:piapprox}

\begin{figure}[htb]
  \centering
  \includegraphics[width=.9\linewidth]{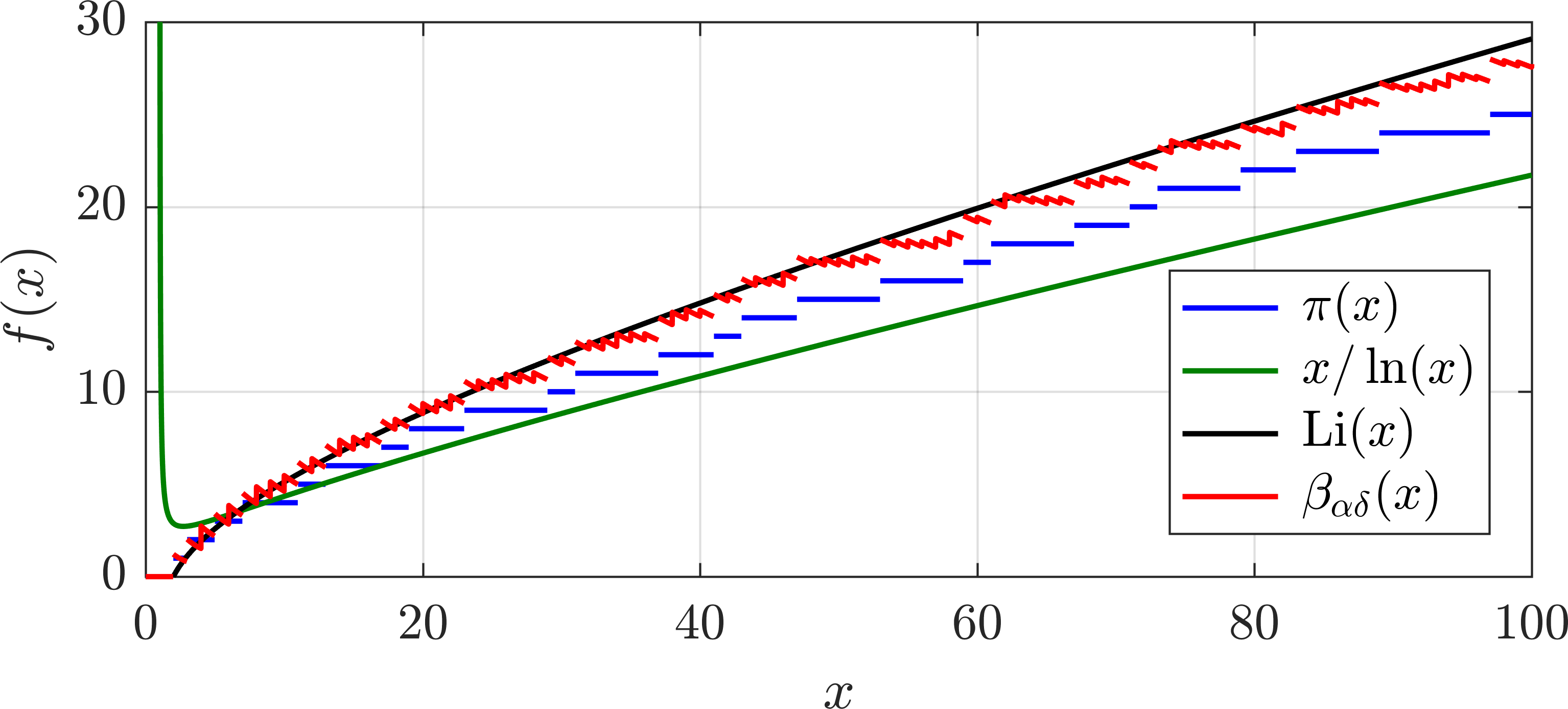}
  \caption{The prime counting function $\pi(x)$ and its approximations.}
  \label{fig:primecount}
\end{figure}

Let $\mathbb{N}:=\{1,2,3,\dots\}$ denote the set of natural numbers,
$\mathbb{N}_0:=\mathbb{N}\cup\{0\}$, and $p_i$ the $i$th prime number with
$i\in\mathbb{N}$, i.e.\ $p_1:=2$, $p_2:=3$, $p_3:=5$, etc.  For
$x\in\mathbb{R}$, the prime counting function is defined as
\begin{equation}
  \label{eq:primecounting}
  \pi(x):=
  \begin{cases}
    0 & \text{if } x<2, \\
    \max \{ k\in \mathbb{N} \,|\, p_k \leq x \} & \text{else}.
  \end{cases}
\end{equation}
According to the prime number theorem, $\pi(x)$ can be approximated by $x/\ln
(x)$, i.e.\ \eqref{eq:piapproxbyln} holds, which will be denoted as
$\pi(x)\sim x/\ln(x)$.  An improved approximation is given by the Eulerian
logarithmic integral \eqref{eq:Li}, i.e.\ $\pi(x)\sim \operatorname{Li}(x)$.
The graphs of the prime counting function and its two approximations are
depicted in Fig.~\ref{fig:primecount} for $x\in[0,100]$.

An alternative approximation of $\pi(x)$ is proposed in
\cite{VartziotisTzavellas2016}. Its definition is based on an additive
function, which will be described briefly in the following. Each $n\in
\mathbb{N}\setminus\{1\}$ can be written as prime factorization $n=\prod_{k\in I(n)}
p_{k}^{\alpha_{k}(n)}$ with the prime numbers $p_k$ and their associated
multiplicities $\alpha_k(n)\in\mathbb{N}$. Here, $I(n)$ denotes the index set
of the prime numbers, which are part of the factorization of $n$. For example,
for $50=2^1\cdot 5^2$ it holds that $I(n)=\{1,3\}$ and $\alpha_1(50)=1$,
$\alpha_3(50)=2$.

An additive function $\beta_{\alpha}:\mathbb{N}\rightarrow \mathbb{N}_0$ is
given by the sum of all prime factors multiplied by their associated
multiplicities, i.e.
\begin{equation}
  \label{eq:betaalpha}
  \beta_{\alpha}(n):=
  \begin{cases}
    0 & \text{for } n=1, \\
    \sum_{k\in I(n)}\alpha_{k}(n)p_{k} & \text{otherwise}.
  \end{cases}
\end{equation}
Here, additive function means that $n_1,n_2\in\mathbb{N}$ implies
$\beta_{\alpha}(n_1n_2)=\beta_{\alpha}(n_1)+\beta_{\alpha}(n_2)$. Furthermore,
for prime numbers it follows readily that $\beta_{\alpha}(p_i)=p_i$,
$i\in\mathbb{N}$.

Summing $\beta_{\alpha}(n)$ for all $n$ less or equal to a given real
number $x>0$ and applying proper scaling leads to the definition of
\begin{equation}
  \label{eq:betaalphadelta}
  \beta_{\alpha\delta}(x):=\frac{12}{\pi^2x}\sum_{n\in\mathbb{N}, n\leq x}\beta_{\alpha}(n)\,.
\end{equation}
The graph of $\beta_{\alpha\delta}(x)$ is depicted red in
Fig.~\ref{fig:primecount}.  As a central result, it has been shown in
\cite{VartziotisTzavellas2016} that $\beta_{\alpha\delta}(x) \sim \pi(x)$,
which is due to the representation
\[
  B_{\alpha}(x):=\sum_{n\in\mathbb{N}, n\leq x}\beta_{\alpha}(n)=
  \frac{\pi^2}{12}\frac{x^2}{\ln x} +{\cal O}\left(\frac{x^2}{\ln^2 x}\right)\,.
\]

This approximation of the prime counting function provides the first
ingredient for deriving prime related fractals. The second ingredient is given
by the following section.

\section{Polygon transformations and Fourier polygons}
\label{sec:polygons}

A geometric transformation scheme for closed polygons based on constructing
similar triangles on the sides is proposed by the authors in
\cite{VartziotisWipper2010LT}. By iteratively applying the transformation,
regular polygons are obtained for specific choices of the transformation
parameters. The proof is based on Fourier polygons, which are also used in
the subsequent section in order to derive the prime number related
fractals. Therefore, this section gives a short introduction to regularizing
polygon transformations and the theoretical background.

\begin{figure}[htbp]
  \centering
  \includegraphics[width=.8\linewidth]{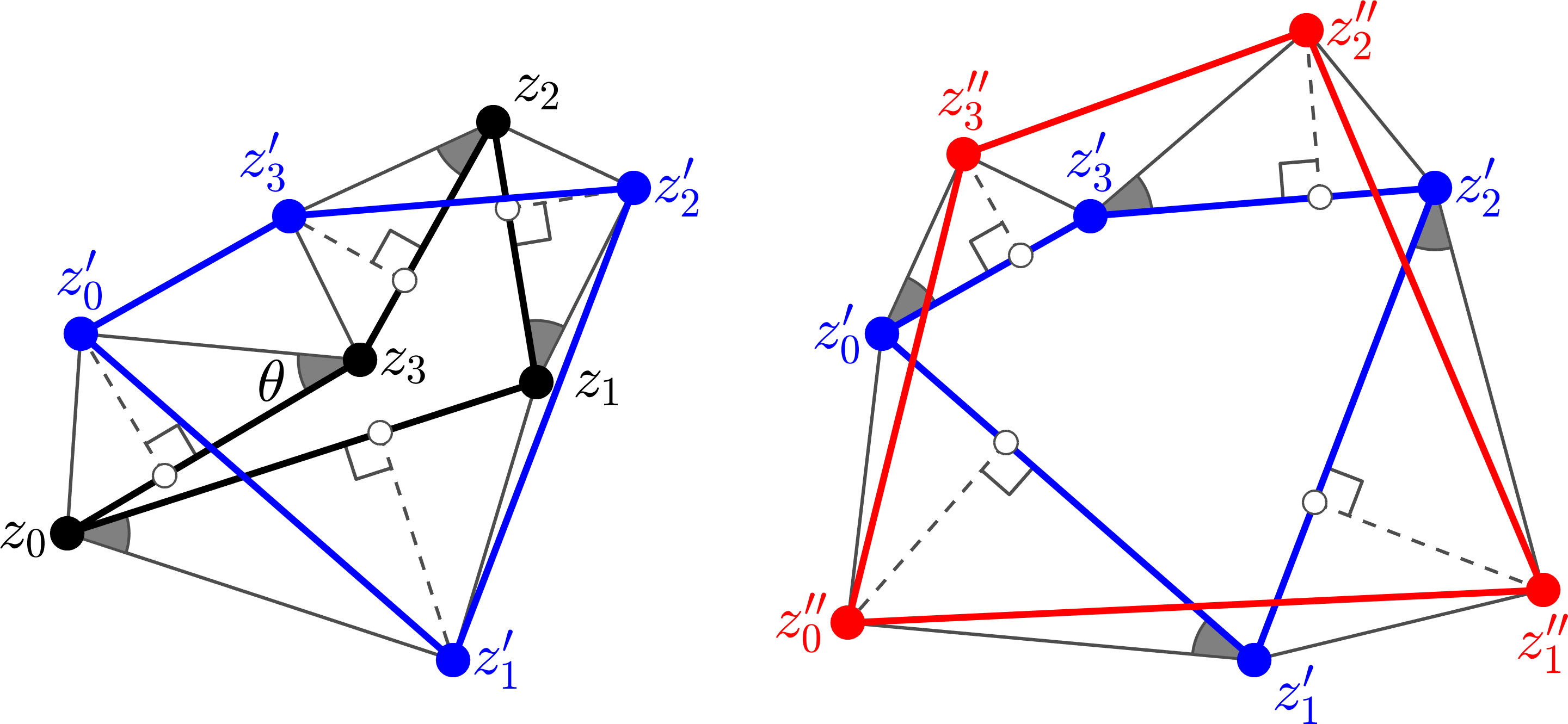}
  \caption{Transformation of an initial polygon by two substeps based on similar triangles.}
  \label{fig:polygontransformation}
\end{figure}

For a given $n\in\mathbb{N}$, $n\geq 3$, let $z=(z_0,\dots,z_{
  n-1})^{\rm t}\in\mathbb{C}^n$ denote an arbitrary polygon in the
complex plane. In the following, all indices have to be taken modulo
$n$. In the first transformation substep, the similar triangles constructed on
each directed side $z_kz_{k+1}$, $k\in\{0,\dots,n-1\}$, of the polygon are
parameterized by a prescribed side subdivision ratio $\lambda\in(0,1)$ and a
base angle $\theta\in(0,\pi/2)$.  This is done by constructing the
perpendicular to the right of the side at the subdivision point $\lambda
z_k+(1-\lambda)z_{k+1}$. On this perpendicular, a new polygon vertex
$z_{k+1}'$ is chosen in such a way that the triangle side $z_kz_{k+1}'$ and
the polygon side $z_kz_{k+1}$ enclose the prescribed angle $\theta$. This is
depicted on the left side of Fig.~\ref{fig:polygontransformation} for an
initial polygon with vertices $z_k$ marked black. The construction of the
first transformation substep leads to a new polygon with vertices $z_k'$
marked blue. The subdivision points, located on the edges of the initial
polygon, are marked by white circles, the associated perpendiculars by dashed
lines, and the angles $\theta$ by grey arcs. In the given example, the
transformation parameters have been set to $\lambda = 1/3$ and $\theta =
\pi/5$.

The rotational effect of the first transformation substep is compensated by
applying a second transformation substep with flipped similar triangles as is
depicted on the right side of Fig.~\ref{fig:polygontransformation}. Starting
from the vertices $z_k'$ of the first transformation substep marked blue, this results in the
polygon with vertices $z_k''$ depicted red. As has been shown in
\cite{VartziotisWipper2010LT}, for
$w:=\lambda+{\rm i}(1-\lambda)\tan\theta$, the new polygon vertices
obtained by applying these two substeps are given by
\[
z_{k+1}':=wz_k+(1-w)z_{k+1}\quad \text{and}\quad
z_{k}'':=(1-\overline{w})z_k'+\overline{w}z_{k+1}'\,,
\]
respectively, where $\overline{w}$ denotes the complex conjugate of $w$.  Both substeps are
linear mappings in $\mathbb{C}^n$ and the combined mapping is given by the
matrix representation
\begin{equation}
  \label{eq:circulanthermitianmatrix}
z''=Mz,\quad\text{where}\quad     M_{j,k}:=
    \begin{cases}
      |1-w|^2+|w|^2 & \text{if } j=k, \\
      w(1-\overline{w}) & \text{if } j=k+1, \\
      \overline{w}(1-w) & \text{if } k=j+1, \\
      0 & \text{otherwise},
    \end{cases}
\end{equation}
with $j,k\in\{0,\dots,n-1\}$.

The transformation matrix $M$ is circulant and Hermitian
\cite{VartziotisWipper2010LT}. Hence, with $r:=\exp(2\pi{\rm i}/n)$
denoting the $n$th root of unity, it holds that $M$ is diagonalized by the
$n\times n$ unitary discrete Fourier matrix
\begin{equation}
  \label{eq:foriermatrix}
F:= \frac{1}{\sqrt{n}}
\begin{pmatrix}
  r^{0\cdot 0} & \dots & r^{0\cdot (n-1)} \cr
  \vdots & \ddots & \vdots \cr
  r^{(n-1)\cdot 0} & \dots & r^{(n-1)\cdot (n-1)} \cr  
\end{pmatrix}
\end{equation}
with entries $F_{j,k}=r^{jk}/\sqrt{n}$ and zero-based indices
$j,k\in\{0,\dots,n-1\}$ \cite{Davis1994}.

Let $z^{(\ell)}:=M^{\ell}z$ denote the polygon obtained by applying the
transformation $\ell$ times. If $\ell$ tends to infinity, the shape of the
scaled limit polygon $z^{(\infty)}$ depends on the dominating eigenvalue of
$M$ and the associated column of $F$. In the following, such a $k$th column of
$F$ is denoted as the $k$th {\em Fourier polygon}, i.e.\
\begin{equation}
  \label{eq:fourierpolygon}
  f_k:=\left(r^{0\cdot k},\dots,r^{(n-1)\cdot k}\right)^{\rm t}, \quad k\in\{0,\dots,n-1\}\,.
\end{equation}
Hence, the Fourier polygons are the prototypes of the limit polygons obtained
by iteratively applying $M$. A full classification of these limit polygons
with respect to the choice of the transformation parameters $\lambda$ and
$\theta$ is given in $\cite{VartziotisWipper2010LT}$. Such transformations
leading to regular polygons can for example be used in finite element mesh
smoothing \cite{VartziotisWipper2009MM}. Furthermore, similar smoothing
schemes can also be applied to volumetric meshes. Here, transformations can
for example be based on geometric constructions
\cite{VartziotisWipper2012Mixed} or on the gradient flow of the mean volume
\cite{VartziotisHimpel2014}.

\begin{figure}[htb]
  \centering
  \includegraphics[width=.9\linewidth]{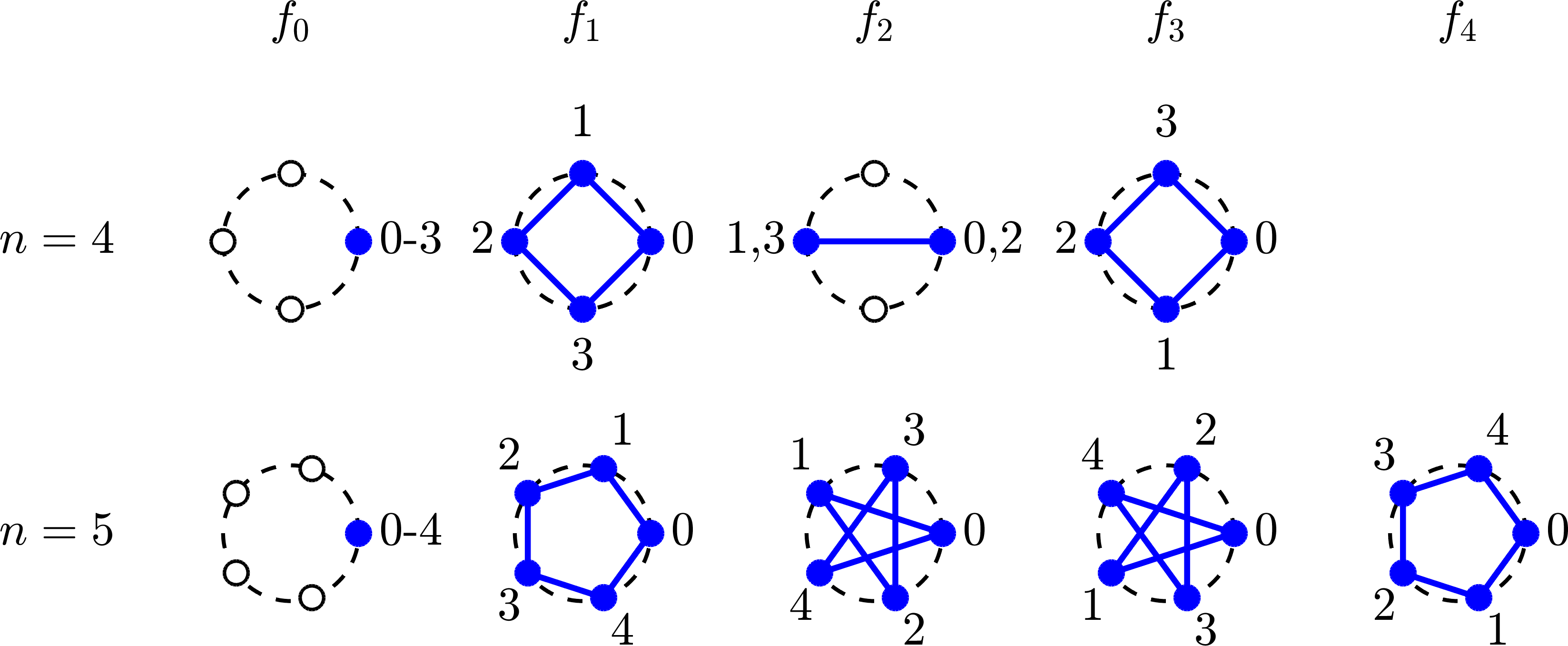}
  \caption{Fourier polygons $f_k$ (blue) for $n=4$ and $n=5$ with given
    vertex indices.}
  \label{fig:fourierpolygons}
\end{figure}

For $n=4$, the four Fourier polygons are depicted in the upper row of
Fig.~\ref{fig:fourierpolygons}. Since the $j$th vertex of $f_k$ is the $jk$th
power of $r$, $f_0$ consists of four times the vertex $1+0{\rm i}$, as is
indicated by the blue point and the vertex label 0-3. The unit circle is marked
by a dashed line, the roots of unity by white circle markers. As can be seen,
$f_1$ is the regular counterclockwise oriented quadrilateral, $f_2$ a segment
with vertex multiplicity two.

In contrast to the case $n=4$, there is no degenerate polygon in the case
$n=5$ except that for $k=0$, as is shown in the lower row of
Fig.~\ref{fig:fourierpolygons}. Since the $k$th Fourier polygon is derived by
iteratively connecting each $k$th vertex defined by the unit roots, and $n=5$ is
a prime number, each polygon vertex has multiplicity one. That is, in the case
of prime numbers, all $f_k$ are either regular $n$-gons or star shaped
$n$-gons if $k>0$. This relation between iterative polygon transformation
limits and prime numbers is also analyzed in
\cite{VartziotisWipper2009PP}. Furthermore, the roots of unity form a cyclic
group under multiplication.

\section{Deriving prime fractals}
\label{sec:fractals}

The graph of the approximate prime counting function $\beta_{\alpha\delta}(x)$
depicted in Fig.~\ref{fig:primecount} gives only an impression for small
values of $x$. Furthermore, it is not suitable to reveal more insight into the
inner structure of prime numbers and their distribution. In search for such a
graphical representation, a combination of the approximate prime counting
function and Fourier polygons is considered in the following.

The main ingredient of the approximate prime counting function $\beta_{\alpha\delta}$
given by \eqref{eq:betaalphadelta} is the sum $B_{\alpha}(x)$. For
$x\in\mathbb{N}\setminus\{1\}$, this sum can be written as
\begin{eqnarray}
  \label{eq:polygonsum}\nonumber
  B_{\alpha}(n)&=&\sum_{k=2}^{n}\beta_{\alpha}(k) = \sum_{k=2}^{n}\sum_{j\in
    I(k)}\alpha_j(k)p_j \\
   &=&\sum_{k\in \mathbb{N}, p_k\leq n}\underbrace{\left(\sum_{j=1}^{\lfloor \log_{p_k}n\rfloor}\left\lfloor \frac{n}{p_k^j}
    \right\rfloor\right)}_{=:w_k(n)}p_k\,,
\end{eqnarray}
where $\lfloor\cdot\rfloor$ denotes rounding towards zero. That is, $B_{\alpha}(n)$ is
a weighted sum of prime numbers. The latter representation can be obtained by
collecting all coefficients in the sum on the left side of
\eqref{eq:polygonsum} for each prime number $p_k$ using a sieve of
Eratosthenes based argument.

The key to prime fractals is to replace the prime numbers $p_k$ in the
representation of $B_{\alpha}(n)$ according to \eqref{eq:polygonsum} by prime number
associated Fourier polygons. This leads to the {\em polygonal prime fractal}
\begin{equation}
  \label{eq:fractalpolygon}
  F_p(n):=\sum_{k\in \mathbb{N}, p_k\leq n} w_k(n)f_{p_k-1}\,,
\end{equation}
with $f_k$ denoting the Fourier polygon according to
\eqref{eq:fourierpolygon}. Here, the Fourier polygon index $p_k-1$ consists of
the $k$th prime number subtracted by one, since zero based indices are used in
the discrete Fourier transformation scheme.

\begin{figure}[htb]
  \centering
  \includegraphics[width=.8\linewidth]{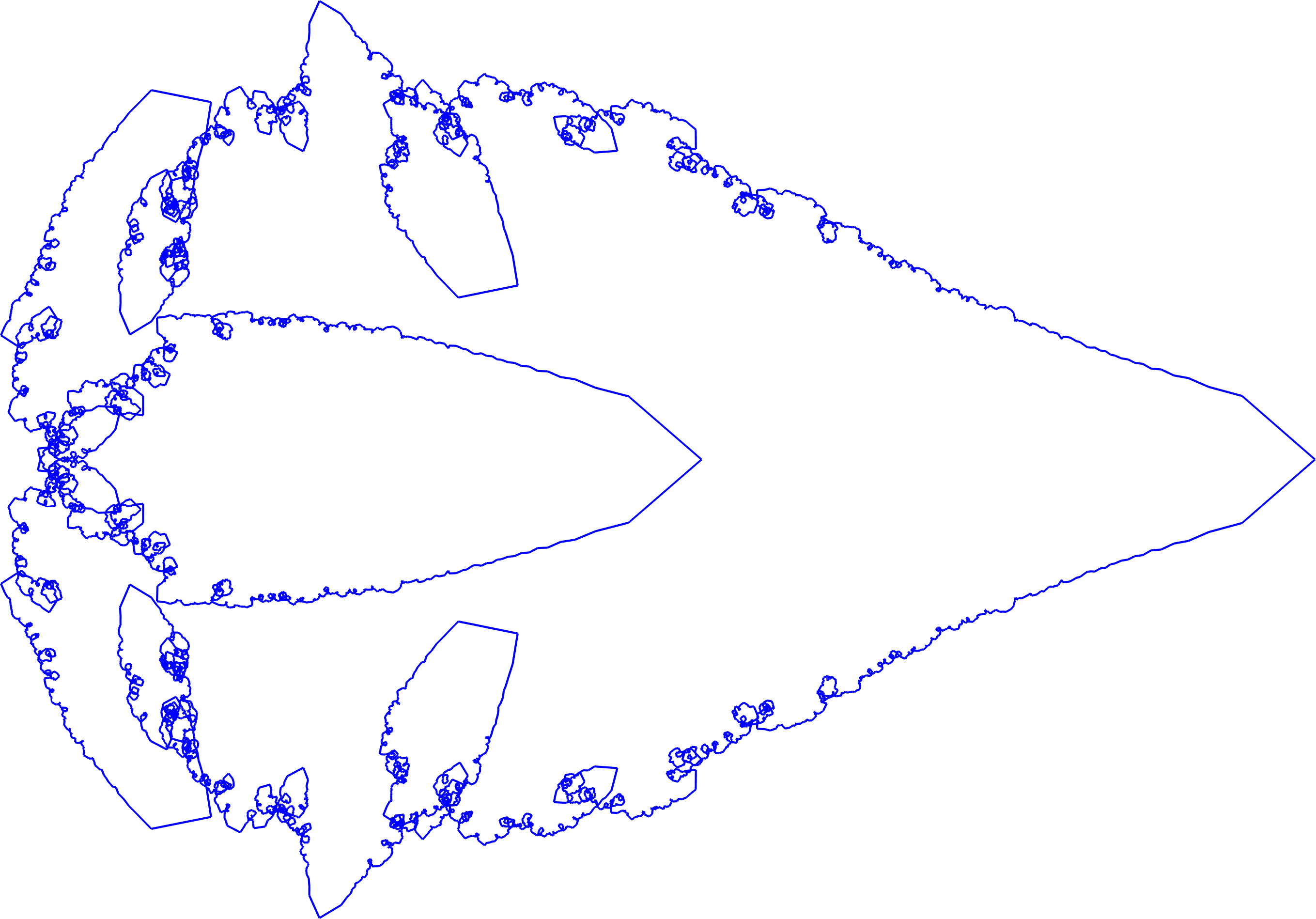}
  \caption{Polygonal prime fractal $F_p(10^4)$.}
  \label{fig:fractalpolygon}
\end{figure}

For $n=10^4$, the polygonal prime fractal $F_p(10^4)$ is depicted in
Fig.~\ref{fig:fractalpolygon}. It is a linear combination of 1\,229 Fourier
polygons $f_{p_k-1}$, each weighted with its associated coefficient $w_k(n)$ as
defined in \eqref{eq:polygonsum}. This polygon consists of $n=10^4$ vertices. The
fractal structure of this polygon is already visible for this comparably low
value of $n$. However, due to its discrete nature, self similarity is not
that obvious for some parts of the polygon. Therefore, an improved fractal is
derived by using the continuous Fourier basis instead of the discrete Fourier
basis. That is, the Fourier polygon $f_k$ is replaced by the Fourier basis
function $f_k(t):=\exp({\rm i}kt)$. This results in the {\em prime fractal
  curve}
\begin{equation}
  \label{eq:fractalcurve}
  F_c(n,t):=\sum_{k\in \mathbb{N}, p_k\leq n} w_k(n)\exp\big({\rm i}(p_k-1)t\big)\,,
  \quad \text{with} \quad t\in(-\pi,\pi]\,.
\end{equation}

\begin{figure}[phtb]
  \centering
  \begin{subfigure}[c]{.95\linewidth}
    \subcaption{Full fractal}
    \includegraphics[width=\linewidth]{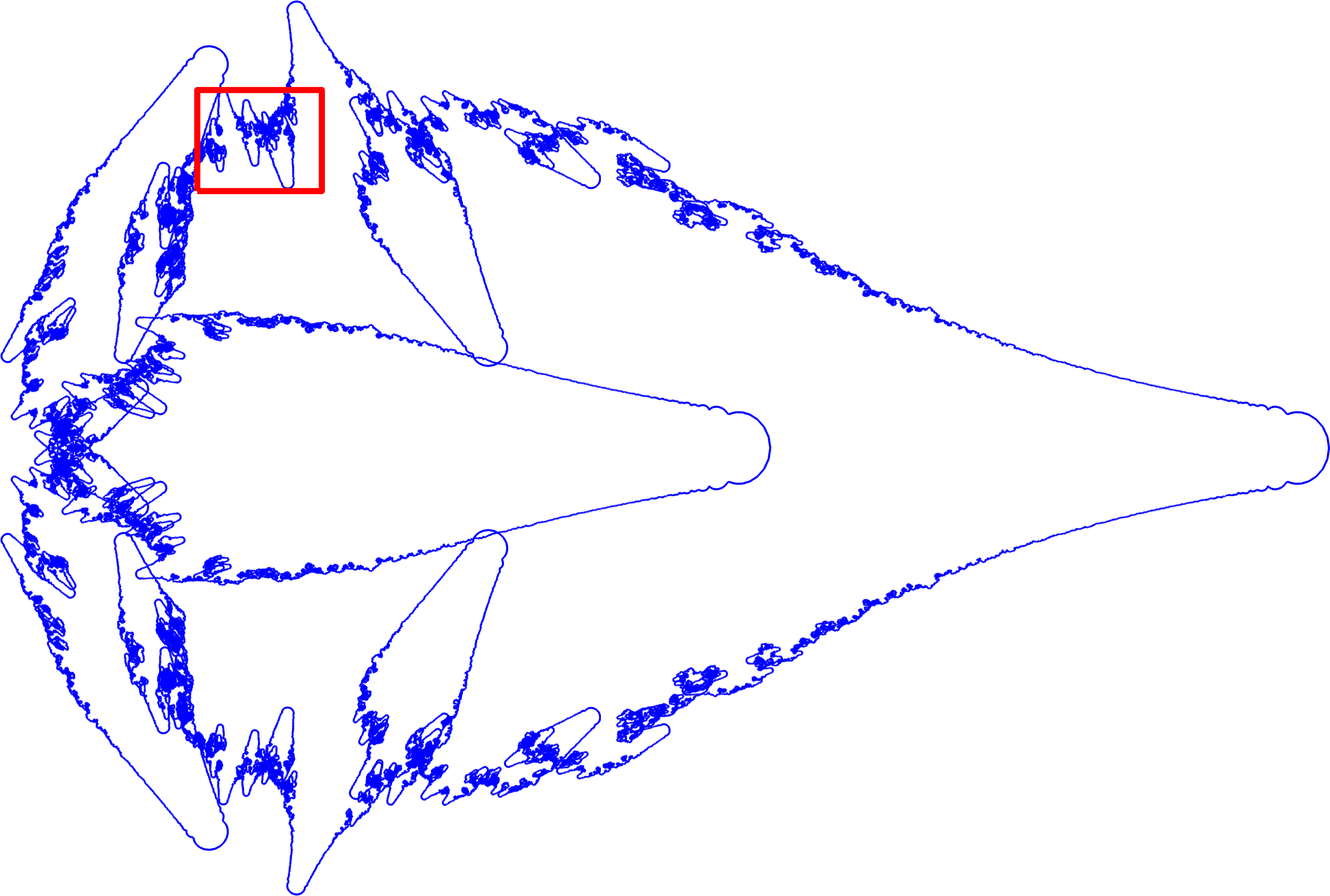}
    \label{fig:fractalcurvefull}
  \end{subfigure}
  \bigskip

  \begin{subfigure}[t]{.48\linewidth}
    \subcaption{Zoom level 1}
    \includegraphics[width=\linewidth]{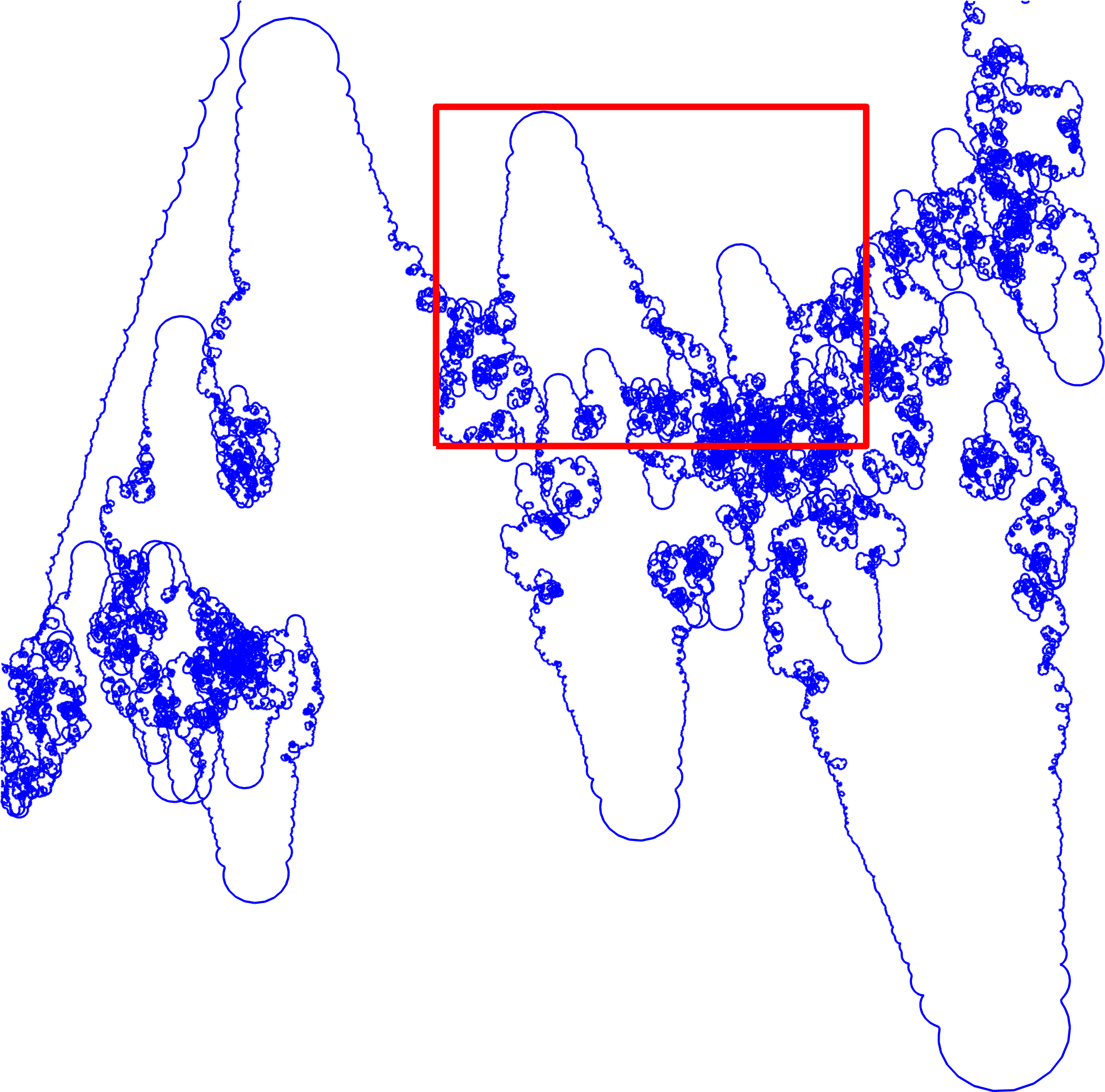}
    \label{fig:fractalcurvezoom1}
  \end{subfigure}
  \hfill
  \begin{subfigure}[t]{.48\linewidth}
    \subcaption{Zoom level 2}
    \includegraphics[width=\linewidth]{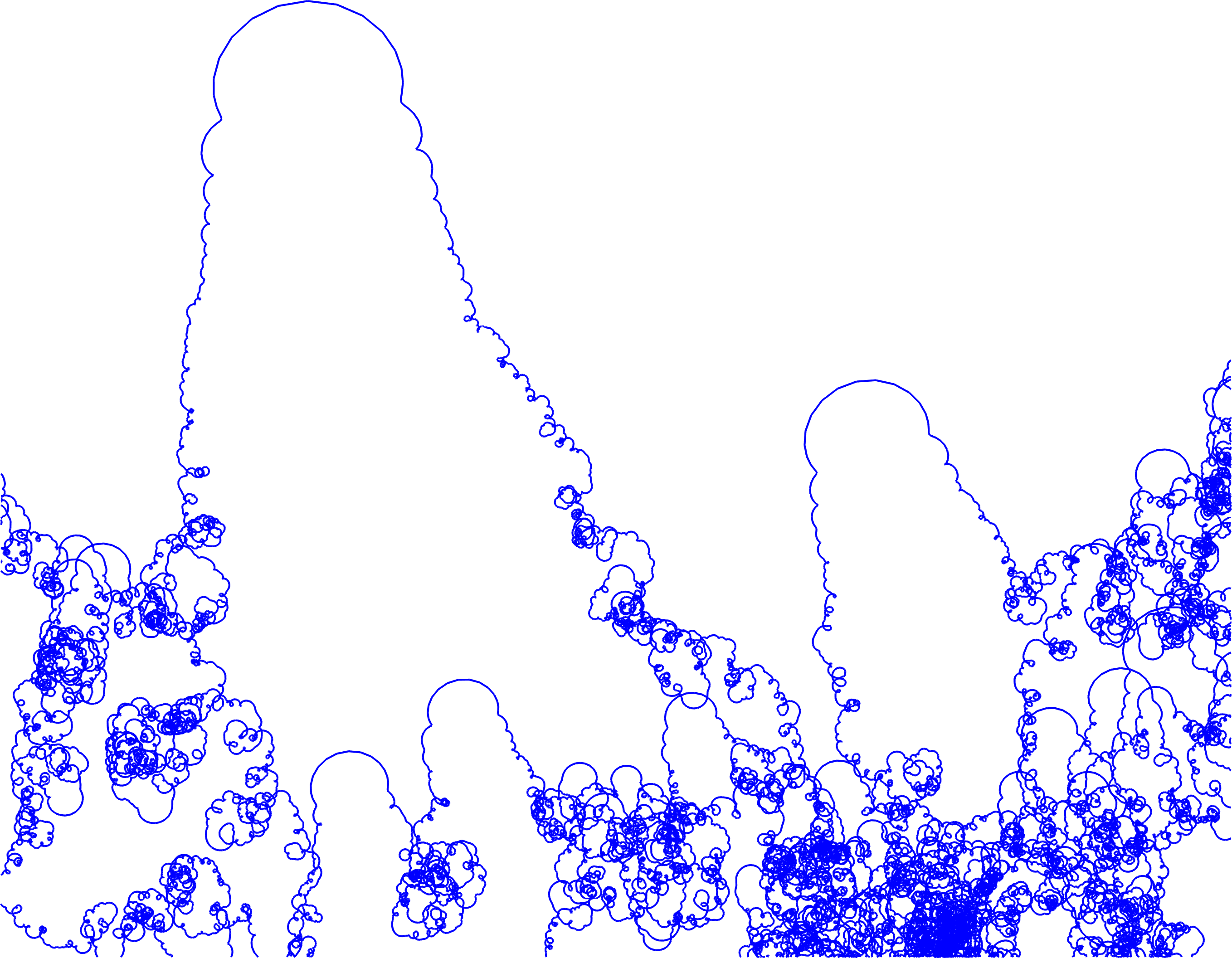}
    \label{fig:fractalcurvezoom2}
  \end{subfigure}

  \caption{Prime fractal curve $F_c$ for $n=10^6$ and $10^7$ evaluation
    points $t_k\in(-\pi,\pi]$. Boxes of successive zooms are marked red.}
\end{figure}
\clearpage

For $n=10^6$, the resulting prime fractal curve is depicted in
Fig.~\ref{fig:fractalcurvefull}. It has been obtained by evaluating $F_c(t)$
at $10^7$ equidistant parameters $t\in(-\pi,\pi]$. In this case, $F_c$ is
the sum of 78\,498 weighted Fourier basis functions. A zoom of the box marked
red in Fig.~\ref{fig:fractalcurvefull} is depicted in
Fig.~\ref{fig:fractalcurvezoom1}. The recurring structures show the fractal
nature of the curve.

\begin{figure}[htb]
  \centering
  \includegraphics[width=.8\linewidth]{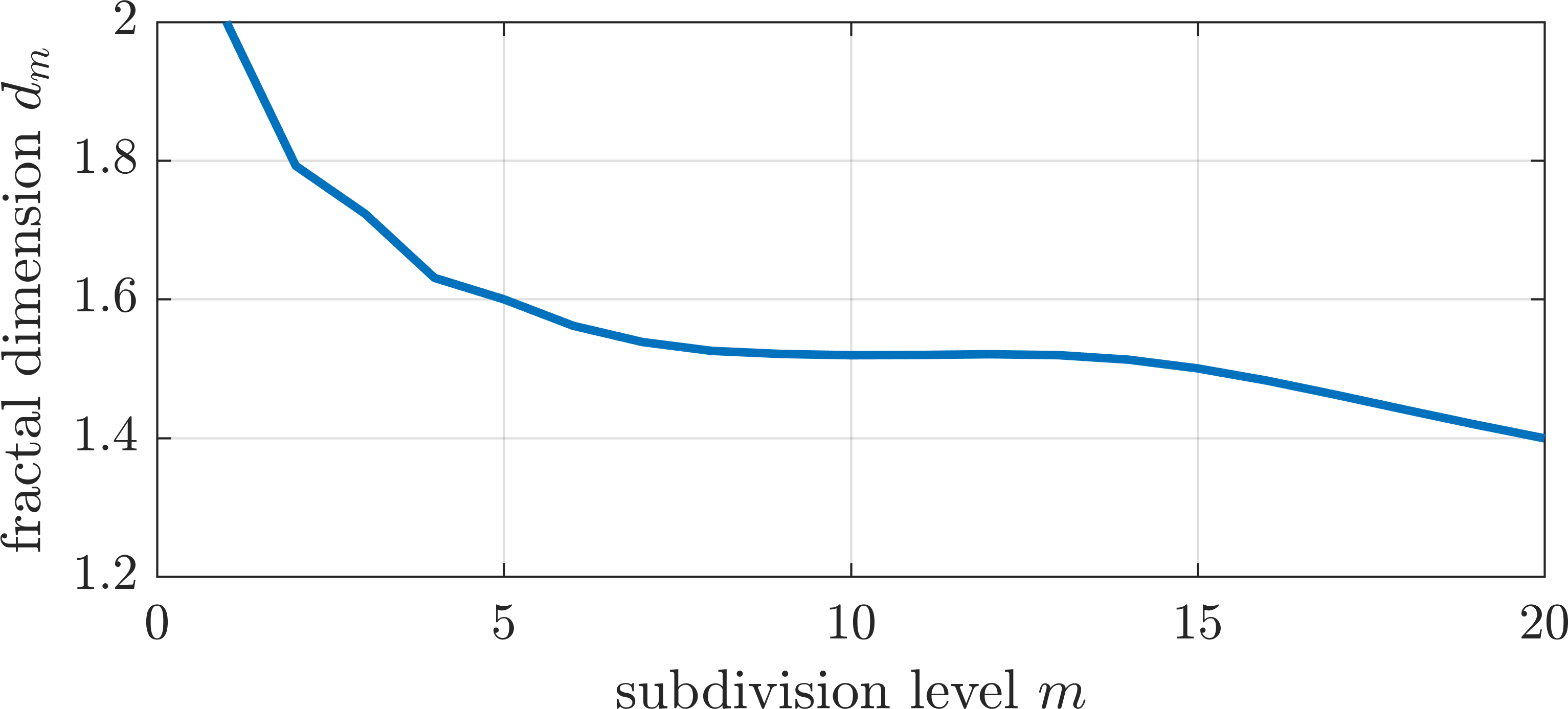}
  \caption{Fractal dimension estimate $d_m$ for the subdivision levels
    $m\in\{0,\dots ,20\}$.}
  \label{fig:fractaldimension}
\end{figure}

The fractal dimension of a curve, also known as Minkowski-Bouligand dimension,
is given by
\begin{equation}
  \label{eq:fractaldimension}
  d:=\lim_{\varepsilon\rightarrow\infty}\frac{\log N(\varepsilon)}{\log(1/\varepsilon)}\,,
\end{equation}
where $\varepsilon$ denotes the edge length of the square boxes covering the
fractal curve and $N(\varepsilon)$ the number of covering boxes
\cite{Falconer2014}. For the curve depicted in
Fig.~\ref{fig:fractalcurvefull}, the fractal dimension is estimated by
recursively subdividing a tight initial square bounding box of $F_c$. For each
subdivision level $m\in\{0,\dots,20\}$, this results in a grid of $2^m\times
2^m$ squares and the associated estimate $d_m:=\log(N_m)/\log(2^m)$, where
$N_m$ denotes the number of grid boxes with nonempty intersections with the
fractal curve. The resulting numerical dimensions are given in
Fig.~\ref{fig:fractaldimension}. Here, the approximate fractal dimension for
the finest subdivision grid results to $d_{20}=1.3995$.

\section{Conclusion}

In this publication, prime number related fractal polygons and curves have
been derived by combining two different results. One is the approximation of
the prime counting function $\pi(x)$ by the partial sum
$\beta_{\alpha\delta}(x)$ based on an additive function as proposed in
\cite{VartziotisTzavellas2016}. The other are prime indexed basis functions of
the continuous Fourier transform. The motivation for this are the column
vectors of the discrete Fourier matrix used for diagonalizing a circulant
Hermitian matrix representing a regularizing transformation of polygons in the
complex plane.

As has been shown in earlier publications \cite{VartziotisWipper2009PP,
  VartziotisWipper2010LT}, there is a relation between the shape of limit
polygons of such iteratively applied polygon transformations and prime
numbers. This is due to the cyclic group defined by the roots of unity and
the exponential representation of the entries of the columns of the discrete
Fourier matrix. The latter are called Fourier polygons. For prime number
indices, these polygons are star shaped.

By replacing the prime numbers in a scaled representation of
$\beta_{\alpha\delta}(x)$ for a given $x\in\mathbb{N}$ by the associated
Fourier polygons, the polygonal prime fractal $F_p$ is derived. Its graphical
features are increased, if $x$ tends to infinity. Alternatively, by using the
prime number related basis functions of the continuous Fourier transform, the
prime fractal curve $F_c$ is derived, which is approximated by $F_p$. The
prime fractal polygon as well as the prime fractal curve show similar patterns
on different scales as has been demonstrated graphically. In addition, a
numerical estimate for the fractal dimension based on the box counting method
was derived for the case $n=10^6$. However, obtaining more detailed results
for much larger $n$ would be desirable.

The given result combines aspects from prime number theory, group theory, and
circulant Hermitian matrices. It is based on the Fourier transformation, which
also plays a role in dynamical systems associated to geometric element
transformations \cite{VartziotisBohnet2016}. The resulting prime fractals
provide an alternative visual representation of an approximation of the prime
counting function and with this of prime numbers and their structures
itself. It is hoped that these alternative representations provide a basis for
further insights into the structure of the distribution of prime
numbers. Furthermore, applying similar visualization techniques to other
number theory functions might reveal additional insights.

\bibliographystyle{hunsrt}
\bibliography{references}
\bigskip

\parindent0cm
Contact: dimitris.vartziotis@nikitec.gr

\end{document}